%&amstex            
\input amstex\documentstyle{amsppt}  
\pagewidth{12.5cm}\pageheight{19cm}\magnification\magstep1
\topmatter
\title Special representations of Weyl groups: a positivity property\endtitle
\author G. Lusztig\endauthor
\address{Department of Mathematics, M.I.T., Cambridge, MA 02139}\endaddress
\thanks{Supported in part by National Science Foundation grant 1303060.}
\endthanks  
\endtopmatter   
\document

\define\Irr{\text{\rm Irr}}

\define\mpb{\medpagebreak}

\define\si{\sim}

\define\sqc{\sqcup}

\define\qua{\quad}

\define\lb{\linebreak}

\define\op{\oplus}
   
\define\part{\partial}
\define\emp{\emptyset}
\define\imp{\implies}

\define\n{\notin}
\define\iy{\infty}
\define\m{\mapsto}
\define\do{\dots}

\define\bsl{\backslash}

\define\sub{\subset}    

\define\T{\times}
\define\ti{\tilde}
\define\nl{\newline}
\redefine\i{^{-1}}

\define\un{\underline}

\define\ot{\otimes}

\define\End{\text{\rm End}}

\define\ind{\text{\rm ind}}

\define\res{\text{\rm res}}

\define\tr{\text{\rm tr}}

\define\he{\heartsuit}

\redefine\c{\chi}
\define\g{\gamma}
\redefine\d{\delta}
\define\e{\epsilon}

\redefine\o{\omega}

\define\ps{\psi}
\define\r{\rho}
\define\s{\sigma}

\define\th{\theta}

\redefine\l{\lambda}

\define\x{\xi}

\redefine\G{\Gamma}

\define\Om{\Omega}

\define\Ps{\Psi}

\define\boc{\bold c}

\define\CC{\bold C}

\define\FF{\bold F}

\define\II{\bold I}
\define\JJ{\bold J}
\define\KK{\bold K}

\define\NN{\bold N}

\define\QQ{\bold Q}
\define\RR{\bold R}

\define\ZZ{\bold Z}

\define\ca{\Cal A}

\define\cd{\Cal D}

\define\ch{\Cal H}

\define\cl{\Cal L}
\define\cm{\Cal M}

\define\fL{\frak L}

\define\fS{\frak S}

\define\BFO{BFO}
\define\EW{EW}
\define\EGNO{EGNO}
\define\KL{KL1}
\define\KLL{KL2}
\define\KM{KM}
\define\SPEC{L1}
\define\ORA{L2}
\define\LEA{L3}
\define\CELLSIV{L4}
\define\POSI{L5}
\define\HEC{L6}
\define\ACTION{L7}
\define\EXCEP{L8}
\define\INV{L9}
\define\PERR{Pe}
\head Introduction\endhead
Let $W$ be an irreducible Weyl group with length function $l:W@>>>\NN$ and let
$S=\{s\in W;l(s)=1\}$. Let $\Irr W$ be a set of representatives for the
isomorphism classes of irreducible representations of $W$ (over $\CC$). In
\cite{\SPEC} a certain subset of $\Irr W$ was defined. The representations in 
this subset were later called {\it special representations}; they play a key 
role in the classification of unipotent representations of a reductive group 
over a finite field $\FF_q$ for which $W$ is the Weyl group. (The definition 
of special representations is reviewed in 3.1.) 

It will be convenient to replace irreducible representations of $W$ with the 
corresponding simple modules of the asymptotic Hecke algebra $\JJ$ (see 
\cite{\HEC, 18.3}) associated to $W$ via the canonical isomorphism 
$\psi:\CC[W]@>\si>>\JJ$ (see 3.1); let $E_\iy$ be the simple $\JJ$-module 
corresponding to $E\in\Irr W$ under $\psi$.
 
In this paper we show that a special representation $E$ of $W$ is characterized
by the following positivity property of $E_\iy$: there exists a $\CC$-basis of
$E_\iy$ such that any element $t_u$ in the standard basis of $\JJ$ acts in 
this basis through a matrix with all entries in $\RR_{\ge0}$. 

The fact that for a special representation $E$, $E_\iy$ has the positivity
property above was pointed out (in the case where $W$ is of classical type) in
\cite{\INV}. In this paper I will recall the argument of \cite{\INV} (see 3.3)
and I give two other proofs which apply for any $W$. One of these proofs (see 
4.4) is based on the interpretation \cite{\LEA}, \cite{\BFO}, of $\JJ$ (or its
part attached to a fixed two-sided cell) in terms of $G$-equivariant vector 
bundles on $X\T X$ where $X$ is a finite set with an action of a finite group 
$G$. Another proof (see Section 2) is based on the use of Perron's theorem for 
matrices with all entries in $\RR_{>0}$. (Previously, Perron's theorem has 
been in used in the context of canonical bases in quantum groups in the study 
\cite{\POSI} of total positivity and, very recently, in the context of the
canonical basis \cite{\KL} of $\CC[W]$, in \cite{\KM}; in both cases the 
positivity properties of the appropriate canonical bases were used). We also 
show that the Hecke algebra representation corresponding to a special
representation $E$ can be realized essentially by a $W$-graph (in the sense of
\cite{\KL}) in which all labels are natural numbers. Some of our results admit
also an extension to the case of affine Weyl groups (see Section 5).

\head 1. Statement of the main theorem\endhead
\subhead 1.1\endsubhead
Let $v$ be an indeterminate and let $\ca=\ZZ[v,v\i]$. 
Let $\ch$ be the Hecke algebra of $W$ that is, the associative $\ca$-algebra 
with $1$ with an $\ca$-basis $\{T_w;w\in W\}$ (where $T_1=1$) and with 
multiplication such that $T_wT_{w'}=T_{ww'}$ if $l(ww')=l(w)+l(w')$ and 
$(T_s+1)(T_s-v^2)=0$ if $s\in S$. Let $\{c_w;w\in W\}$ be the $\ca$-basis of 
$\ch$ denoted by $\{C'_w;w\in W\}$ in \cite{\KL} (with $q=v^2$); see also
\cite{\HEC, 5.2}. For example, if $s\in S$, we have $c_s=v\i T_s+v\i$. The 
left cells and two-sided cells of $W$ are the equivalence classes for the 
relations $\si_L$ and $\si_{LR}$ on $W$ defined in \cite{\KL}, see also
\cite{\HEC, 8.1}; we shall write $\si$ instead of $\si_L$. For $x,y$ in $W$ we
have $c_xc_y=\sum_{z\in W}h_{x,y,z}c_z$ where $h_{x,y,z}\in\NN[v,v\i]$. As in 
\cite{\HEC, 13.6}, for $z\in W$ we define $a(z)\in\NN$ by
$h_{x,y,z}\in v^{a(z)}\ZZ[v\i]$ for all $x,y$ in $W$ and
$h_{x,y,z}\n v^{a(z)-1}\ZZ[v\i]$ for some $x,y$ in $W$. (For example, 
$a(1)=0$ and $a(s)=1$ if $s\in S$.) For $x,y,z$ in $W$ we have 
$h_{x,y,z}=\g_{x,y,z\i}v^{a(z)}\mod v^{a(z)-1}\ZZ[v\i]$ where 
$\g_{x,y,z\i}\in\NN$ is well defined.
Let $\JJ$ be the $\CC$-vector space with basis $\{t_w;w\in W\}$. For $x,y$ in 
$W$ we set $t_xt_y=\sum_{z\in W}\g_{x,y,z\i}t_z\in\JJ$. This defines a 
structure of associative $\CC$-algebra on $\JJ$ with unit element of the form 
$\sum_{d\in\cd}t_d$ where $\cd$ is a certain subset of the set of involutions 
in $W$, see \cite{\HEC, 18.3}. 
For any subset $X$ of $W$ let $\JJ_X$ be the subspace of $\JJ$ with basis 
$\{t_w;w\in X\}$; let $\JJ^+_X$ be the set of elements of the form 
$\sum_{w\in X}f_wt_w\in\JJ_X$ with $f_w\in\RR_{>0}$ for all $w\in X$. We have 
$\JJ=\op_\boc\JJ_\boc$ where $\boc$ runs over the two-sided cells of $W$. Each
$\JJ_\boc$ is a subalgebra of $\JJ$ with unit element $\sum_{d\in\cd_\boc}t_d$
where $\cd_\boc=\boc\cap\cd$; moreover, $\JJ_\boc\JJ_{\boc'}=0$ if 
$\boc\ne\boc'$. 

{\it Until the end of Section 4 we fix a two-sided cell $\boc$.} 
\nl
Let $L$ be the set of left cells that are contained in $\boc$. We have 
$$\boc=\sqc_{\G\in L}\G=\sqc_{\G,\G'\text{ in }L}(\G\cap\G'{}\i);$$
moreover, $\G\cap\G'{}\i\ne\emp$ for any $\G,\G'\text{ in }L$. It follows that
$$\JJ_\boc=\op_{\G\in L}\JJ_\G=\op_{\G,\G'\text{ in }L}\JJ_{\G\cap\G'{}\i};$$
moreover, $\JJ_{\G\cap\G'{}\i}\ne0$. Note that for $\G\in L$, $\JJ_\G$ is a
left ideal of $\JJ_\boc$. 

A line $\cl$ in $\JJ_{\G\cap\G'{}\i}$ is said to be {\it positive} if 
$\cl^+:=\cl\cap\JJ^+_{\G\cap\G'{}\i}\ne\emp$; in this case, $\cl^+$ consists 
of all $\RR_{>0}$-multiples of a single nonzero vector. We now state our main 
result.

\proclaim{Theorem 1.2}(a) Let $\G\in L$. There is a unique left ideal $M_\G$ 
of $\JJ_\boc$ such that property ($\he$) below holds:

($\he$) $M_\G=\op_{\G'\in L}M_{\G,\G'}$ where for any $\G'\in L$,  
$M_{\G,\G'}:=M_\G\cap\JJ_{\G\cap\G'{}\i}$ is a positive line.

(b) Let $\G\in L,\G'\in L,u\in\boc$. We have $u\in\ti\G\cap\ti\G'{}\i$ for 
well-defined $\ti\G,\ti\G'$ in $L$. If $\ti\G\ne\G'$, then $t_uM_{\G,\G'}=0$. 
If $\ti\G=\G'$, then $t_uM_{\G,\G'}=M_{\G,\ti\G'}$ and 
$t_uM_{\G,\G'}^+=M_{\G,\ti\G'}^+$.

(c) The $\JJ_\boc$-module $M_\G$ in (a) is simple. Its isomorphism class is 
independent of $\G\in L$.

(d) The subspace $\II=\op_{\G\in L}M_\G$ of $\JJ_\boc$ is a simple two-sided
ideal of $\JJ_\boc$.

(e) Let $\G,\G',\ti\G,\ti\G'$ be in $L$. If $\G\ne\ti\G'$ then
$M_{\G,\G'}M_{\ti\G,\ti\G'}=0$. If $\G=\ti\G'$, then multiplication in 
$\JJ_\boc$ defines an isomorphism 
$M_{\G,\G'}\ot M_{\ti\G,\G}@>\si>>M_{\ti\G,\G'}$ and a surjective map 
$M_{\G,\G'}^+\T M_{\ti\G,\G}^+@>>>M_{\ti\G,\G'}^+$.

(f) Let $\G,\G'$ be in $L$. The antiautomorphism $\th:\JJ_\boc@>>>\JJ_\boc$ 
given by $t_x\m t_{x\i}$ for all $x\in\boc$ maps $M_{\G,\G'}$ onto 
$M_{\G',\G}$ and $M_{\G,\G'}^+$ onto $M_{\G',\G}^+$.
\endproclaim
The proof is given in Section 2.

\subhead 1.3\endsubhead
As a consequence of Theorem 1.2, the simple $\JJ_\boc$-module $M_\G$ admits a 
$\CC$-basis $\{\ti e_{\G'};\G'\in L\}$ with the following property:

(i) {\it If $u\in\boc$ and $\G'\in L$, then $t_u\ti e_{\G'}$ is an
$\RR_{\ge0}$-linear combination of elements $\ti e_{\G''}$ with $\G''\in L$; 
more precisely, if $u\in\ti\G\cap\ti\G'{}\i$ with $\ti\G,\ti\G'$ in $L$, then

$t_u\ti e_{\G'}=\l_{u,\G',\ti\G'}\ti e_{\ti\G'}$
\nl
with $\l_{u,\G',\ti\G'}\in\RR_{>0}$ if $\ti\G=\G'$ and $\l_{u,\G',\ti\G'}=0$ 
if $\ti\G\ne\G'$.}
\nl
Indeed, we can take for $\ti e_{\G'}$ any element of $M_{\G,\G'}^+$ and we use 
1.2(b).

\subhead 1.4\endsubhead
Let $\le$ be the standard partial order on $W$. By \cite{\KL}, to any $y\ne w$
in $W$ one can attach a number $\mu(y,w)\in\ZZ$ such that for any $s\in S$ and
any $w\in W$ with $sw>w$ we have $c_sc_w=\sum_{y\in W;sy<y}\mu(y,w)c_y$. By 
\cite{\KLL} we have $\mu(y,w)\in\NN$.

\subhead 1.5\endsubhead
Let $\un\ch=\CC(v)\ot_\ca\ch$ where we use the obvious imbedding
$\ca@>>>\CC(v)$; we denote $1\ot c_w$ again by $c_w$. Let 
$\un\JJ=\CC(v)\ot_\CC\JJ$ where we use the obvious imbedding
$\CC@>>>\CC(v)$. We have a homomorphism of $\CC(v)$-algebras (with $1$)
$\Ps:\un\ch@>>>\un\JJ$ given by 
$$\Ps(c_x)=\sum_{d\in\cd,z\in W,d\si z}h_{x,d,z}t_z$$
for all $x\in W$, see \cite{\HEC, 18.9}. (Note that $\Ps$ is in fact the
composition of a homomorphism in {\it loc.cit.} with an automorphism of 
$\un\ch$.)

\subhead 1.6\endsubhead
For any $\G\in L$ let $S_\G$ be the set of all $t\in S$ such that $rt<r$ for 
some (or equivalently any) $r\in\G$.

We fix $\G\in L$. Let $\{\ti e_{\G'};\G'\in L\}$ be a $\CC$-basis of $M_\G$ as
in 1.3; we use the notation of 1.3. We shall view $\CC(v)\ot M_\G$ as an 
$\un\ch$-module via $\Ps$. Let $s\in S$ and let $\G'\in L$; let $\d$ be the 
unique element in $\G'\cap\cd$. We show:

(a) {\it If $s\in S_{\G'}$, then $\Ps(T_s)\ti e_{\G'}=v^2\ti e_{\G'}$.}

(b) {\it If $s\n S_{\G'}$, then 
$$\Ps(T_s)\ti e_{\G'}=
-\ti e_{\G'}+\sum_{\ti\G\in L;s\in S_{\ti\G}}f_{\ti\G,\G'}v\i\ti e_{\ti\G}$$
where}
$$f_{\ti\G,\G'}=\sum_{u\in\G'\cap\ti\G\i}\mu(u,\d)\l_{u,\G',\ti\G}\in
\RR_{\ge0}.$$
\nl
By definition, we have
$$\Ps(c_s)\ti e_{\G'}=\sum_{d\in\cd,u\in c,d\si u}h_{s,d,u}t_u\ti e_{\G'}=
\sum_{\ti\G\in L}\sum_{d\in\cd,u\in\G'\cap\ti\G\i,d\in\G'}
h_{s,d,u}\l_{u,\G',\ti\G}\ti e_{\ti\G}.$$
Since in the last sum we have $d\in\G'$ we see that we can assume that $d=\d$.
Thus we have
$$\Ps(c_s)\ti e_{\G'}=\sum_{\ti\G\in L}\sum_{u\in\G'\cap\ti\G\i}
h_{s,\d,u}\l_{u,\G',\ti\G}\ti e_{\ti\G}.$$
If $s\d<\d$ (that is, $s\in S_{\G'}$) we have $c_sc_\d=(v+v\i)c_\d$ hence 
$h_{s,\d,u}$ is $(v+v\i)$ for $u=\d$ and is $0$ for $u\ne\d$; hence in this 
case
$$\Ps(c_s)\ti e_{\G'}=(v+v\i)\ti e_{\G'};$$
(we use that $\l_{\d,\G',\G'}=1$.)

We now assume that $s\d>\d$ (that is, $s\n S_{\G'}$). In this case, 
$h_{s,\d,u}$ is $\mu_{u,\d}$ if $su<u$ and is $0$ if $su>u$ (see 1.4); hence 
$$\align&\Ps(c_s)\ti e_{\G'}=\sum_{\ti\G\in L}\sum_{u\in\G'\cap\ti\G\i;su<u}
\mu(u,\d)\l_{u,\G',\ti\G}\ti e_{\ti\G}\\&=\sum_{\ti\G\in L;s\in S_{\ti\G}}
\sum_{u\in\G'\cap\ti\G\i}\mu(u,\d)\l_{u,\G',\ti\G}\ti e_{\ti\G}.\endalign$$
Now (a),(b) follow. 

Note that (a),(b) show that in the $\un\ch$-module $\CC(v)\ot M_\G$ the 
generators $T_s$ act with respect to the basis $\{\ti e_{\G'};\G'\in L\}$ 
essentially by formulas which are those in a $W$-graph (in the sense of 
\cite{\KL}) in which all labels are in $\RR_{\ge0}$.

\subhead 1.7\endsubhead
In Section 4 we will give another proof of the existence part of 1.2(a) which 
also shows that $\ti e_{\G'}$ in 1.3 can be chosen so that 

(i) each $\ti e_{\G'}$ is a $\ZZ_{>0}$-linear combination of elements in 
$\{t_x;x\in\G\cap\G'{}\i\}$,

(ii) $\l_{u,\G',\G'_1}\in\ZZ_{\ge0}$ (notation of 1.3).
\nl
In particular, with this choice of $\ti e_{\G'}$, the constants 
$f_{\ti\G,\G'}$ in the "$W$-graph formulas" in 1.6 are in $\ZZ_{\ge0}$.

\head 2. Proof of Theorem 1.2\endhead
\subhead 2.1\endsubhead
From \cite{\HEC,\S15} we see that, for $x,y,u$ in $W$ we have:
$$\g_{x,y,u}=\g_{y,u,x}=\g_{u,x,y},\tag a$$
$$\g_{x,y,u}\ne0\imp x\si y\i,y\si u\i,u\si x\i.\tag b$$
By \cite{\HEC, 18.4(a)}:

(c) {\it for $y,z$ in $W$ we have $y\si z$ if and only if $t_yt_{z\i}\ne0$.}

\subhead 2.2\endsubhead
Let $\G,\G',\ti\G,\ti\G'$ be in $L$. From 2.1(b) we deduce:
$$\text{ If }\G\ne\ti\G',\text{ then }
\JJ_{\G\cap\G'{}\i}\JJ_{\ti\G\cap\ti\G'{}\i}=0,\tag a$$
$$\JJ_{\G\cap\G'{}\i}\JJ_{\ti\G\cap\G\i}\sub\JJ_{\ti\G\cap\G'{}\i}.\tag b$$
We show:
$$\text{If }u\in\G\cap\G'{}\i,\text{ then }
t_u\JJ_{\ti\G\cap\G\i}^+\sub\JJ_{\ti\G\cap\G'{}\i}^+.\tag c$$
Let $\x=\sum_{y\in\ti\G\cap\G\i}f_yt_y\in\JJ_{\ti\G\cap\G\i}$ with 
$f_y\in\RR_{>0}$ for all $y$. We must show that 
$t_u\x\in\JJ_{\ti\G\cap\G'{}\i}^+$; it is enough to show that for any 
$z\in\ti\G\cap\G'{}\i$ there exists $y\in\ti\G\cap\G\i$ such that
$\g_{u,y,z\i}\ne0$ or that there exists $y\in W$ such that $\g_{z\i,u,y}\ne0$ 
(see 2.1(a)); such $y$ is automatically in $\ti\G\cap\G\i$. Hence it is enough
to show that for any $z\in\ti\G\cap\G'{}\i$ we have $t_{z\i}t_u\ne0$. This 
holds since $z\i\si u\i$ (see 2.1(c)).

\mpb

From (c) we deduce
$$\JJ_{\G\cap\G'{}\i}^+\JJ_{\ti\G\cap\G\i}^+\sub\JJ_{\ti\G\cap\G'{}\i}^+.
\tag d$$

\subhead 2.3\endsubhead
Let $\G\in L$. For any $\G'\in L$ we define a $\CC$-linear map 
$T_{\G'}:\JJ_{\G\cap\G'{}\i}@>>>\JJ_{\G\cap\G'{}\i}$ by
$$T_{\G'}(t_x)=\sum_{y\in\G\cap\G\i}t_xt_y=
\sum_{y\in\G\cap\G\i,z\in\G}\g_{x,y,z\i}t_z.$$
We show:

(a) {\it the matrix representing $T_{\G'}$ with respect to the basis
$\{t_w;w\in\G\cap\G'{}\i\}$ has all entries in $\RR_{>0}$.}
\nl
An equivalent statement is: for any $x,z$ in $\G\cap\G'{}\i$, the sum
$\sum_{y\in\G\cap\G\i}\g_{x,y,z\i}$ is $>0$. Since $\g_{x,y,z\i}\in\NN$ for
all $y$, it is enough to show that for some $y\in\G\cap\G\i$ we have 
$\g_{x,y,z\i}\ne0$ or equivalently (see 2.1(a)) that for some $y\in W$ we have
$\g_{z\i,x,y}\ne0$ (we then have automatically $y\in\G\cap\G\i$). Thus, it is 
enough to show that $t_{z\i}t_x\ne0$. This follows from 2.1(c) since 
$z\i\si x\i$.

\mpb

Applying Perron's theorem \cite{\PERR} to the matrix in (a) we see that there 
is a unique $T_{\G'}$-stable positive line $\cl_{\G,\G'}$ in 
$\JJ_{\G\cap\G'{}\i}$ (the "Perron line").

Now let $u\in\boc$; we have, $u\in\ti\G\cap\ti\G'{}\i$ with $\ti\G,\ti\G'$ in 
$L$. From 2.2(a), 2.2(d), we deduce 

(b) If $\ti\G\ne\G'$, then $t_u\JJ_{\G\cap\G'{}\i}=0$ hence 
$t_u\cl_{\G,\G'}=0$;

(c) if $\ti\G=\G'$, then
$t_u\cl_{\G,\G'}^+\sub\JJ^+_{\G\cap\ti\G'{}\i}$, hence $t_u\cl_{\G,\G'}$ 
is a positive line in $\JJ_{\G\cap\ti\G'{}\i}$.
\nl
In the setup of (c), we have $t_u(T_{\G'}(\x))=T_{\ti\G'}(t_u\x)$ for any 
$\x\in\JJ_{\G\cap\G'{}\i}$. It follows that $t_u\cl_{\G,\G'}$ is a 
$T_{\ti\G'}$-stable line in $\JJ_{\G\cap\ti\G'{}\i}$. Thus,

(d) if $\ti\G=\G'$, then $t_u\cl_{\G,\G'}=\cl_{\G,\ti\G'}$.
\nl
We set $\cm_\G=\op_{\G'\in L}\cl_{\G,\G'}$. From (b),(d) we see that $\cm_\G$
is a $\JJ_\boc$-submodule of $\JJ_\G$.

We now see that the existence part of 1.2(a) is proved: we can take 
$M_\G=\cm_\G$.

\subhead 2.4\endsubhead
Let $\G\in L$ and let $M_\G$ be any $\JJ_\boc$-submodule of $\JJ_\G$ for which
property ($\he$) in 1.2(a) holds. We show that 1.2(b) holds for $M_\G$. Let 
$u\in\ti\G\cap\ti\G'{}\i$ be as in 1.2(b) and let $\G'\in L$. If 
$\ti\G\ne\G'$, then $t_u\JJ_{\G\cap\G'{}\i}=0$ hence $t_uM_{\G,\G'}=0$. Now 
assume that $\ti\G=\G'$. By 2.2(c), left multiplication by $t_u$ maps 
$\JJ^+_{\G\cap\G'{}\i}$ into $\JJ^+_{\G\cap\ti\G'{}\i}$ hence it maps any
positive line in $\JJ_{\G\cap\G'{}\i}$ onto a positive line in 
$\JJ_{\G\cap\ti\G'{}\i}$. In particular, it maps $M_{\G,\G'}$ onto a line in 
$\JJ_{\G\cap\ti\G'{}\i}$, which, being also contained in $M_\G$, must be equal
to $M_{\G,\ti\G'}$; moreover, it maps $M_{\G,\G'}^+$ into
$\JJ^+_{\G\cap\ti\G'{}\i}$ hence onto $M_{\G,\ti\G'}^+$. Thus, 1.2(b) holds for
$M_\G$.

We now choose a basis $\{\ti e_{\G'};\G'\in L\}$ of $M_\G$ such that
$\ti e_{\G'}\in M_{\G,\G'}^+$ for any $\G'\in L$; then for any $u\in\boc$, the
matrix of the $t_u$-action on $M_\G$ in this basis has entries in $\RR_{\ge0}$.
Thus,

(a) $\tr(t_u,M_\G)\in\RR_{\ge0}$ for all $u\in\boc$.
\nl
We show:

(b){\it the $\CC$-linear map $\nu:\JJ_\boc@>>>\End_\CC(M_\G)$ given by the 
$\JJ_\boc$-module structure on $M_\G$ is surjective.}
\nl
It is enough to show that for any $\G',\ti\G'$ in $L$ there exists $u\in\boc$ 
such that $\nu(t_u)$ carries the line $M_{\G,\G'}$ onto the line 
$M_{\G,\ti\G'}$ and carries the line $M_{\G,\G''}$ (where $\G''\in L$,
$\G''\ne\G'$) to zero. Note that any $u\in\G'\cap\ti\G'{}\i$ has the required 
properties. This proves (b).

It follows that the $\JJ_\boc$-module $M_\G$ is simple. We show:

(c) {\it Assume that $M'$ is any simple $\JJ_\boc$-module such that
$\tr(t_u,M')\in\RR_{\ge0}$ for all $u\in\boc$. Then $M'$ is isomorphic to 
$M_\G$.}
\nl
Assume that this is not so. We use the orthogonality formula 
$$\sum_{u\in\boc}\tr(t_u,M_\G)\tr(t_{u\i},M')=0,$$
which is a special case of \cite{\HEC, 19.2(e)} (taking into account 
\cite{\HEC, 20.1(b)} and using that $u\m u\i$ maps $\boc$ into $\boc$). Since 
each term in the last sum is in $\RR_{\ge0}$, it follows that each term in the
last sum is $0$. In particular, we have $\tr(t_d,M_\G)\tr(t_d,M')=0$ for any 
$d\in\cd_\boc$. We show that for any $d\in\cd_\boc$ we have 
$\tr(t_d,M_\G)\in\RR_{>0}$. Using the basis of $M_\G$ employed in the proof of
(a), it is enough to show that some diagonal entry of the matrix of the
$t_d$-action in this basis is $\ne0$ (all entries are in $\RR_{\ge0}$). We 
have $d\in\G'$ for a unique $\G'\in L$; then $t_dM_{\G,\G'}^+=M_{\G,\G'}^+$ 
and the desired property holds.

From $\tr(t_d,M_\G)\tr(t_d,M')=0$ and $\tr(t_d,M_\G)\in\RR_{>0}$ we deduce 
that $\tr(t_d,M')=0$ for any $d\in\cd_\boc$. Since $\sum_{d\in\cd_\boc}t_d$ is 
the unit element $1_\boc$ of $\JJ_\boc$, it follows that $\tr(1_\boc,M')=0$. 
This is a contradiction. This proves (c).

\mpb

Let $I$ be the simple ideal of $\JJ_\boc$ such that $I\cm_\G\ne0$. It is a
$\CC$-vector space of dimension $N^2$ where $N$ is the number of elements in 
$L$. If $\ti\G\in L$, then $\cm_{\ti\G}$ is a simple $\JJ_\boc$-module such 
that $\tr(t_u,\cm_{\ti\G})\in\RR_{\ge0}$ for all $u\in c$ (we use (a) with
$M_\G$ replaced by $\cm_{\ti\G}$); hence, by (c), we have 
$\cm_{\ti\G}\cong M_\G$ as $\JJ_\boc$-modules. In particular, the isomorphism 
class of $\cm_{\ti\G}$ is independent of $\ti\G$. We see that the (necessarily
direct) sum $\sum_{\ti\G\in L}\cm_{\ti\G}$ is contained in $I$ and has 
dimension $N^2$ hence it is equal to $I$; we also see that $M_\G\sub I$ and, 
taking intersections with $\JJ_\G$, we see that $M_\G\sub\cm_\G$, hence 
$M_\G=\cm_\G$ (since $\dim M_\G=\dim\cm_\G=N$). We now see that the uniqueness
part of 1.2(a) is proved. Note that 1.2(b),1.2(c),1.2(d) are also proved and 
we have $\II=I$.

\subhead 2.5\endsubhead
We prove 1.2(e). In the setup of (e), if $\G\ne\ti\G'$ then, using 2.2(a), we
have $M_{\G,\G'}M_{\ti\G,\ti\G'}=0$. Assume now that $\G=\ti\G'$. Using 2.2(d),
we see that $M_{\G,\G'}^+M_{\ti\G,\G}^+$ is contained in
$\JJ_{\ti\G\cap\G'{}\i}^+$; it is also contained in $\II$ (since $\II$ is 
closed under multiplication), hence it is contained in 
$\II\cap\JJ_{\ti\G\cap\G'{}\i}^+=M_{\ti\G,\G'}^+$. Thus, multiplication 
restricts to a map $M_{\G,\G'}^+\T M_{\ti\G,\G}^+@>>>M_{\ti\G,\G'}^+$. This 
map is necessarily surjective since $M_{\ti\G,\G'}^+$ is a single orbit of 
$\RR_{>0}$ under scalar multiplication. This implies that the linear map
between lines $M_{\G,\G'}\ot M_{\ti\G,\G}@>>>M_{\ti\G,\G'}$ is an isomorphism.
This proves 1.2(e).

\subhead 2.6\endsubhead
We prove 1.2(f). Any element $\x\in\JJ_\boc$ defines a linear map 
${}^t(\th(\x)):M_\G^*@>>>M_\G^*$ where $M_\G^*$ denotes the dual space and 
${}^t$ denotes the transpose. This defines a $\JJ_\boc$-module structure on 
$M_\G^*$ such that for any $x\in\boc$ we have 
$\tr(t_x,M_\G^*)=\tr(t_{x\i},M_\G)$. By the argument in \cite{\HEC, 20.13(a)},
$\tr(t_{x\i},M_\G)$ is the complex conjugate of $\tr(t_x,M_\G)$. But the last 
trace is a real number, so that $\tr(t_x,M_\G^*)=\tr(t_x,M_\G)$. It follows
that $M_\G^*\cong M_\G$ as $\JJ_\boc$-modules. From the definitions, the 
simple two-sided ideal $\II'$ of $\JJ_\boc$ such that $\II'M_\G^*\ne0$  
satisfies $\II'=\th(\II)$. It follows that $\th(\II)=\II$. Since 
$\II=\op_{\ti\G,\ti\G'\text{ in }L}M_{\ti\G,\ti\G'}$ and
$\th(\JJ_{\G,\G'})=\JJ_{\G',\G}$, it follows that
$$\th(M_{\G,\G'})\sub
\JJ_{\G',\G}\cap\op_{\ti\G,\ti\G'\text{ in }L}M_{\ti\G,\ti\G'}=M_{\G',\G}.$$
Since $\th$ is a vector space isomorphism, it follows that
$\th(M_{\G,\G'})=M_{\G',\G}$. Note that $\th(\JJ_{\G,\G'}^+)=\JJ_{\G',\G}^+$; 
hence 
$$\th(M_{\G,\G'}^+)\sub\JJ_{\G',\G}^+\cap M_{\G',\G}=M_{\G',\G}^+.$$
This forces the equality $\th(M_{\G,\G'}^+)=M_{\G',\G}^+$ (since 
$M_{\G',\G}^+$ is a single orbit of $\RR_{>0}$ under scalar multiplication). 
This proves 1.2(f). Theorem 1.2 is proved.

\subhead 2.7\endsubhead
After an earlier version of this paper was posted, P. Etingof told me that the
line $M_{\G,\G'}$ in 1.2 is the same as the line associated in 
\cite{\EGNO, 3.4.4} to the right $\JJ_{\G\cap\G\i}$-module 
$\JJ_{\G\cap\G'{}\i}$ (viewed as a based module over a based ring) that is, 
the unique positive line $\fL$ in $\JJ_{\G\cap\G'{}\i}$ such that $\fL$ is a 
right $\JJ_{\G\cap\G\i}$-submodule. (The discussion in {\it loc.cit.} concerns
left (instead of right) indecomposable based modules over a fusion ring.) 
Indeed, from the definitions we see that $\fL$ must be the same as 
$\cl_{\G,\G'}$ in 2.3, hence the same as $M_{\G,\G'}$.

\head 3. Special representations\endhead
\subhead 3.1\endsubhead
When $\ch$ is tensored with $\CC$ (using the ring homomorphism $\ca@>>>\CC$,
$v\m1$), then it becomes $\CC[W]$, the group algebra of $W$. (For $w\in W$ we 
have $1\ot T_w=w\in\CC[W]$; we denote $1\ot c_w$ again by $c_w$.) We have a 
homomorphism of $\CC$-algebras (with $1$) $\ps:\CC[W]@>>>\JJ$ given by 
$$\ps(c_x)=\sum_{d\in\cd,z\in W,d\si z}h_{x,d,z}|_{v=1}t_z$$
for all $x\in W$, see \cite{\HEC, 18.9}; this is an isomorphism, see
\cite{\HEC, 20.1}. For example, if $W=\{1,s\}$ is of type $A_1$ we have
$\ps(c_1)=t_1+t_s$, $\ps(c_s)=2t_s$; hence $\ps(1)=t_1+t_s$, 
$\ps(s)=-t_1+t_s$.

For each $E\in\Irr W$ let $E_\iy$ be the simple $\JJ$-module corresponding
to $E$ under $\ps$ and let $\boc_E$ be the unique two-sided cell of $W$ such 
that $\JJ_{\boc_E}E\ne0$. (Note that $E=E_\iy$ as $\CC$-vector spaces.) Let 
$\Irr^\boc W=\{E\in\Irr W;\boc_E=\boc\}$ and let $a'=a(xw_0)$ for any 
$x\in\boc$, where $w_0$ is the longest element of $W$.

For any $k\in\NN$ let $\fS^k$ be the $k$-th symmetric power of the reflection
representation of $W$, viewed as a representation of $W$ in an obvious way.
For $E\in\Irr W$ let $b_E$ be the smallest integer $k\ge0$ such that $E$ is a
constituent of $\fS^k$.
Now for any $E\in\Irr^\boc W$ we have $b_E\ge a'$ and there is a unique 
$E\in\Irr^\boc W$ such that $b_E=a'$; this $E$ is denoted by $E^\boc$ and is 
called the {\it special representation} associated to $\boc$. (This is a 
reformulation of the definition of special representations given in 
\cite{\SPEC}.) 

\proclaim{Theorem 3.2}In the setup of Theorem 1.2, for any $\G\in L$, we have
$M_\G\cong E^\boc_\iy$ as $\JJ_\boc$-modules.
\endproclaim
We give two proofs; one is contained in 3.3,3.4,3.5. The other is given in
3.4,3.6.

\subhead 3.3\endsubhead
In this subsection we assume that $W$ is of type $A,B$ or $D$.

Let $\G\in L$. For any $\G'\in L$ we set
$$\e_{\G'}=\sum_{z\in\G\cap\G'{}\i}t_z\in\JJ^+_{\G\cap\G'{}\i}.$$
By \cite{\INV, 4.8(b)}, $\{\e_{\G'};\G'\in L\}$ is a $\CC$-basis of the 
unique $\JJ$-submodule of $\JJ_\G$ isomorphic to $E^\boc_\iy$. 
By the uniqueness part of 1.2(a) this $\JJ$-submodule of $\JJ_\G$ (viewed
as a $\JJ_\boc$-module) must be the same as $M_\G$ in 1.2(a). We see that in 
this case, $M_\G$ is isomorphic to $E^\boc_\iy$ and $M_{\G,\G'}$ is the line 
spanned by $\e_{\G'}$. In particular, 3.2 holds in our case.

\subhead 3.4\endsubhead
In this subsection we assume that $\boc$ is such that $\Irr^\boc W$ consists 
of exactly $2$ irreducible representations. In this case, $W$ is of type $E_7$ 
(resp. $E_8$) and the $2$ irreducible representations in $\Irr^\boc W$ have 
degree $512$ (resp. $4096$). Let $\G\in L$ and let $d\in\cd\cap\G$. The 
$\CC$-linear map $r:\JJ_\G@>>>\JJ_\G$ given by left multiplication by 
$(-1)^{l(d)}\psi(w_0)$ is in fact $\JJ$-linear (since $w_0$ is central in $W$)
and $r(t_x)=t_{x^*}$ for any $x\in\G$, where $x\m x^*$ is a certain fixed 
point free involution of $\G$, see \cite{\ACTION}. Then
$\JJ_\G^1=\{\x\in\JJ_\G;r(\x)=\x\}$ is a simple $\JJ_\boc$-submodule of 
$\JJ_\G$ with $\CC$-basis $\{t_x+t_{x^*};x\in\G_1\}$ where $\G_1$ is a set of 
representatives for the orbits of $x\m x^*$ on $\G$. Note that, if $x\in\G$,
then $\{x,x^*\}$ is the intersection of $\G$ with the inverse of a left cell 
$\G'\in L$; hence $t_x+t_{x^*}\in\JJ^+_{\G\cap\G'{}\i}$. By the uniqueness 
part of 1.2(a), we must have $M_\G=\JJ_\G^1$ and that for any $\G'\in L$, 
$M_{\G,\G'}$ is the line spanned by $\sum_{x\in\G\cap\G'{}\i}t_x$.

Now let $E\in\Irr^\boc W$ be such that $E_\iy=\JJ_\G^1$. From the definitions
we have $\tr((-1)^{l(d)}w_0,E)=|\G_1|=\dim E$. Hence, if $\e=\pm1$ is the 
scalar by which $w_0$ acts on $E$, then $(-1)^{l(d)}\e=1$. We have
$l(d)=a(d)\mod2$; hence $\e=(-1)^{a(d)}$. But this equality characterizes
the special representation in $\Irr^\boc W$ (the special representation 
satisfies it, the nonspecial representation doesn't satisfy it). We see that 
$M_\G=\JJ_\G^1\cong E^\boc_\iy$. In particular, 3.2 holds in our case.

\subhead 3.5\endsubhead
In this subsection we assume that $W$ is of exceptional type, but that $\boc$ 
is not as in 3.4. In this case, $E^\boc$ is the only representation in 
$\Irr^\boc W$ of dimension equal to $|L|$; since $M_\G$ (in 1.2(a)) has 
dimension equal to $|L|$, it follows that $M_\G=\JJ_\G^1\cong E^c_\iy$. In 
particular, 3.2 holds in our case. This completes the proof of Theorem 3.2.

\subhead 3.6\endsubhead
In this subsection we give a second proof of Theorem 3.2 assuming that $\boc$ 
is not as in 3.4. Let $a=a(w)$ for any $w\in\boc$. Let $\G\in L$. Let 
$X=\sum_{w\in W}v^{-l(w)}T_w\in\un\ch$. We can view $\CC(v)\ot\JJ_\G$ and 
$\CC(v)\ot M_\G$ as $\un\ch$-modules via $\Ps$ in 1.5. By \cite{\INV, 4.6}, 
for any $x\in\G$ we have
$$Xt_x=v^a\sum_{z\in\boc}t_zt_x\mod\sum_{i<a}v^i\JJ_\G.$$
By an argument as in the proof of 2.2(c), we see that
$\sum_{z\in\boc}t_zt_x\in\JJ_\G^+$. It follows that if $\G'\in L$ and
$\x\in M_{\G,\G'}$ then 
$$X\x=v^a\x'\mod\sum_{i<a}v^i\JJ_\G$$
where $\x'\in \JJ_\G^+$. In particular, we have $X\x\ne0$. Thus,
$X(\CC(v)\ot M_\G)\ne0$. Using this and Theorem 4.2 in \cite{\INV} we deduce
that the simple $\un\ch$-module $\CC(v)\ot M_\G$ is a constituent of the 
"involution module" $M$ in \cite{\INV, 0.1} (with $\QQ(u)$ replaced by
$\CC(v)$). According to \cite{\EXCEP} if a simple $\un\ch$-module appears in 
$\CC(v)\ot\JJ_{\G'}$ for every $\G'\in L$ and it appears in $M$, then that 
$\un\ch$-module corresponds to $E^\boc$. We deduce that $M_\G\cong E^\boc_\iy$.
This completes the second proof of Theorem 3.2, assuming that $\boc$ is not as
in 3.4.

\head 4. Equivariant vector bundles\endhead
\subhead 4.1\endsubhead
In this section we fix a reductive, not necessarily connected algebraic group 
$G$ over $\CC$ acting on a finite set $X$. Let $G\bsl X$ be the set of 
$G$-orbits on $X$. Representations of reductive groups over $\CC$ are always 
assumed to be of finite dimension over $\CC$ and algebraic. For $x\in X$ let 
$G_x=\{g\in G;gx=x\}$. Now $G$ acts 
diagonally on $X\T X$ and we can consider the Grothendieck group $K_G(X\T X)$ 
of $G$-equivariant complex vector bundles ($G$-eq.v.b.) on $X\T X$ . This is 
an (associative) ring with $1$ under convolution, denoted by $*$ (see 
\cite{\LEA, 2.2}, \cite{\CELLSIV, 10.2}). For a $G$-eq.v.b. $V$ on $X\T X$ we 
denote by $V_{x,y}$ the fibre of $V$ at $(x,y)\in X\T X$. Let $B$ be the set 
of pairs $(\Om,\r)$ where $\Om$ is a $G$-orbit $\Om$ in $X\T X$ and $\r$ is an 
irreducible representation of $G_\Om$ (the isotropy group of a point 
$(x,y)\in\Om$). For any $(\Om,\r)\in B$ we denote by $V^{\Om,\r}$ the 
$G$-eq.v.b. on $X\T X$ such that $V_{\Om,\r}|_{X\T X-\Om}=0$ and the action of
$G_\Om$ on $V^{\Om,\r}_{x,y}$ is equivalent to $\r$. Now 
$\{V^{\Om,\r};(\Om,\r)\in B\}$ is a $\ZZ$-basis of $K_G(X\T X)$. Let 
$\KK_G(X\T X)=\CC\ot K_G(X\T X)$, viewed as a $\CC$-algebra.

\subhead 4.2\endsubhead
In this subsection we assume that $G$ is finite. Let $\o,\o'$ be in $G\bsl X$.
Let $V^{\o,\o'}$ be the $G$-eq.v.b. on $X\T X$ such that 
$V^{\o,\o'}_{a,b}=\CC[G]$ if $(a,b)\in\o\T\o'$, $V^{\o,\o'}_{a,b}=0$ if 
$(a,b)\n\o\T\o'$. (Here $\CC[G]$ is the left regular representation of $G$.) 
The $G$-action $g:V^{\o,\o'}_{a,b}@>>>V^{\o,\o'}_{ga,gb}$ is left translation 
by $g$ on $\CC[G]$ (if $(a,b)\in\o\T\o'$) and is $0$ if $(a,b)\n\o\T\o'$. We 
show:

(a) {\it Let $(Om,\r)\in B$; we have $\Om\sub\o_1\T\o'_1$ where $\o_1,\o'_1$ 
are in $G\bsl X$. Then $U':=V^{\Om,\r}*V^{\o,\o'}$ is isomorphic to a direct 
sum of copies of the single $G$-eq.v.b. $V^{\o_1,\o'}$. More precisely, if 
$\o'_1\ne\o$, we have $U'=0$; if $\o'_1=\o$, we have 
$U'=V_{\o_1,\o'}^{\op(\dim\r|\Om||\o_1|\i)}$.}
\nl
For $(a,b)\in X\T X$ we have
$$U'_{a,b}=\op_{z\in\o;(a,z)\in\Om}V^{\Om,\r}_{a,z}\ot\CC[G]\text{ if }
b\in\o',$$
$$U'_{a,b}=0\text{ if }b\n\o'.$$
Thus the support of $U'$ is contained in $\o_1\T\o'$ and $U'=0$ unless 
$\o'_1=\o$. We now assume that $\o'_1=\o$ and $(a,b)\in\o_1\T\o'$. Then
$U'_{a,b}=\op_{z;(a,z)\in\Om}V^{\Om,\r}_{a,z}\ot\CC[G]$. We have
$\dim U'_{a,b}=d|G||\Om|/|\o_1|$. We show:

(b) {\it as a $G_a\cap G_b$-module, $U'_{a,b}$ is a multiple of the regular 
representation.}
\nl 
Let $\s_1,\do,\s_k$ be the various $G_a\cap G_b$-orbits contained in $\o'$. We 
have $U'_{a,b}=\op_{i=1}^kR_i$, where 
$R_i=\op_{z\in\s_i}V^{\Om,\r}_{a,z}\ot\CC[G]$. We pick $z_i\in\s_i$. Now 
$$R_i=\ind_{G_a\cap G_b\cap G_{z_i}}^{G_a\cap G_b}(A\ot B)$$
where 
$$A=\res_{G_a\cap G_b\cap G_{z_i}}^{G_a\cap G_{z_i}}(V^{\Om,\r}_{a,z_i}),\qua
B=\res_{G_a\cap G_b\cap G_{z_i}}^{G_{z_i}\cap G_b}(\CC[G]).$$
It is enough to show that $R_i$ is a multiple of the regular representation of
$G_a\cap G_b$. Since $R_i$ is induced, it is enough to show that $A\ot B$ is a
multiple of the regular representation of $G_a\cap G_b\cap G_{z_i}$. It is 
also enough to show that $B$ is a multiple of the regular representation of 
$G_a\cap G_b\cap G_{z_i}$. This follows from the fact that $\CC[G]$ is a 
multiple of the regular representation of $G_{z_i}\cap G_b$. This proves (b). 
Now (a) follows.

\mpb

Note that in $\KK_G(X\T X)$ we have
$$V^{\o,\o'}=\sum_{(\Om,\r)\in B;\Om\sub\o\T\o'}\dim\r|\Om|V^{\Om,\r}.\tag c$$

\subhead 4.3\endsubhead
We now drop the assumption (in 4.2) that $G$ is finite. Let $\hat\KK_G(X\T X)$
be the $\CC$-vector space consisting of formal (possibly infinite) linear 
combinations $\sum_{(\Om,\r)\in B}f_{\Om,\r}V^{\Om\r}$ where
$f_{\Om,\r}\in\CC$. The left $\KK_G(X\T X)$-module structure on $\KK_G(X\T X)$
given by left multiplication extends naturally to a left $\KK_G(X\T X)$-module
structure on $\hat\KK_G(X\T X)$. If $\o,\o'$ are as in $G\bsl X$ then we can 
define $V^{\o,\o'}\in\hat\KK_G(X\T X)$ by the sum 4.2(c) (which is now a 
possibly infinite sum); we set $\bar V^{\o,\o'}=|\o|\i V^{\o,\o'}$ so that
$$\bar V^{\o,\o'}=
\sum_{(\Om,\r)\in B;\Om\sub\o\T\o'}\dim\r|\Om||\o|\i V^{\Om,\r}.\tag a$$
Now formula 4.2(a) extends to the present case as follows. Let $(Om,\r)\in B$;
we have $\Om\sub\o_1\T\o'_1$ where $\o_1,\o'_1$ are $G$-orbits in $X$. Then
$$V^{\Om,\r}V^{\o,\o'}=NV^{\o_1,\o'}\tag b$$
where $N\in\ZZ$ is $0$ if $\o'_1\ne\o$ and $N=\dim\r|\Om||\o_1|\i$ if
$\o'_1=\o$. Hence
$$V^{\Om,\r}\bar V^{\o,\o'}=N'\bar V^{\o_1,\o'}\tag c$$
where $N'\in\ZZ$ is $0$ if $\o'_1\ne\o$ and $N'=\dim\r|\Om||\o|\i$ if
$\o'_1=\o$. 

\mpb

Let $\o'\in G\bsl X$. Let $\ti R_{\o'}$ be the subspace of $\hat\KK_H(X\T X)$ 
consisting of formal (possibly infinite) linear combinations
$\sum_{(\Om,\r)\in B;pr_2\Om=\o'}f_{\Om,\r}V^{\Om,\r}$ with $f_{\Om,\r}\in\CC$.
Here $pr_2:X\T X@>>>X$ is the second projection. Note that $\ti R_{\o'}$ is a 
$\KK_G(X\T X)$-submodule of $\hat\KK_H(X\T X)$.

Let $R_{\o'}$ be the subspace of $\hat\KK_G(X\T X)$ with basis formed by the 
elements $\bar V^{\o,\o'}$ for various $\o\in G\bsl X$. Using (c) we see that 
$R_{\o'}$ is a (simple) $\KK_G(X\T X)$-submodule of $\hat\KK_G(X\T X)$; we 
have $R_{\o'}\sub\ti R_{\o'}$. Using (c) we see also that if
$\o''\in G\bsl X$ then $\bar V^{\o,\o'}\m\bar V^{\o,\o''}$ defines an
isomorphism of $\KK_G(X\T X)$-modules $R_{\o'}@>\si>>R_{\o''}$. Hence the 
isomorphism class of the $\KK_G(X\T X)$-module $R_{\o'}$ is independent of the
choice of $\o'$.

\subhead 4.4\endsubhead
We now assume that $G$ is the finite group associated to $\boc$ in \cite{\ORA}
and that $X$ is the finite $G$-set $\op_{\G\in L}G/H_\G$ where $H_\G$ is the 
subgroup of $G$ defined in \cite{\LEA, 3.8}. 
In this case we have $\hat\KK_G(X\T X)=\KK_G(X\T X)$. By a conjecture in
\cite{\LEA, 3.15}, proved in \cite{\BFO}, there exists an isomorphism of 
$\CC$-algebras $\c:\KK_G(X\T X)@>\si>>\JJ_\boc$ carrying the basis 
$(V^{\Om,\r})$ of $\KK_G(X\T X)$ onto the basis $\{t_x;x\in\boc\}$ of 
$\JJ_\boc$. Under $\c$, the left ideal $\JJ_\G$ of $\JJ_\boc$ (for $\G\in L$) 
corresponds to the left ideal $\ti R_{\o'}$ of $\KK_G(X\T X)$ where 
$\o'\in G\bsl X$ corresponds to $\G$, and the basis $\{t_x;x\in\G\}$ of 
$\JJ_\G$ corresponds to the intersection of the basis $(V^{\Om,\r})$ of 
$\KK_G(X\T X)$ with $\ti R_{\o'}$. The basis of $R_{\o'}$ formed by the 
elements $\bar V^{\o,\o'}$ corresponds to a family of elements 
$\{e_{\G'};\G'\in L\}$ in $\JJ_\G$.

From 4.3(a) we see that $e_{\G'}\in\JJ^+_{\G\cap\G'{}\i}$ for any $\G'\in L$
(in fact the coefficients of the various $t_x, x\in\G\cap\G'{}\i$ are in
$\ZZ_{>0})$ and from 4.3(c) we see that for $u\in\boc$ the product
$t_ue_{\G'}$ is a $\ZZ_{\ge0}$ multiple of an element $e_{\G''}$. We see that 
the $\CC$-subspace of $\JJ_\G$ spanned by $\{e_{\G'};\G'\in L\}$ satisfies 
property ($\he$) in 1.2(a) hence it is equal to $M_\G$. This provides another 
proof in our case for the existence part of 
1.2(a), with the additional integrality properties in 1.7.

\head 5. Final remarks\endhead
\subhead 5.1\endsubhead
Theorem 1.2 and its proof remain valid if $W$ is replaced by an affine Weyl
group (with $\boc$ assumed to be finite) or by a finite Coxeter group; in the
last case we use the positivity property of $h_{x,y,z}$ established in 
\cite{\EW}. In these cases, the simple $\JJ_\boc$-module given by Theorem 1.2 
will be called the special $\JJ_\boc$-module.

\subhead 5.2\endsubhead
Assume now that $W$ is an (irreducible) affine Weyl group and that $\boc$ is a
not necessarily finite two-sided cell of $W$. We denote again by $L$ the set 
of left cells of $W$ that are contained in $\boc$; this is a finite set. Then 
the $\CC$-algebra $\JJ_\boc$ with its basis $\{t_x;x\in\boc\}$ is defined. Let
$\hat JJ_\boc$ be the set of formal (possibly infinite) linear combinations
$\sum_{u\in\boc}f_ut_u$ where $f_u\in\CC$. This is naturally a left
$\JJ_\boc$-module. For any subset $X$ of $\boc$ let $\hat\JJ_X$ be the set of 
all $\sum_{u\in\boc}f_ut_u\in\hat\JJ_\boc$ such that $f_u=0$ for $u\in\boc-X$.
If $\G\in L$, then $\hat\JJ_\G$ is a $\JJ_\boc$-submodule of $\hat\JJ_\boc$. 
We have $\hat\JJ_\G=\op_{\G'\in L}\hat\JJ_{\G\cap\G'{}\i}$.

According to a conjecture in \cite{\CELLSIV, 10.5}, proved in \cite{\BFO}, we
can find $G,X$ as in 4.1 and an isomorphism of $\CC$-algebras
$\c:\KK_G(X\T X)@>\si>>\JJ_\boc$ carrying the basis 
$\{V^{\Om,\r};(\Om,\r)\in B\}$ of $\KK_G(X\T X)$ onto the basis 
$\{t_x;x\in\boc\}$ of $\JJ_\boc$. This extends in an obvious way to an 
isomorphism $\hat\c:\hat\KK_G(X\T X)@>\si>>\hat\JJ_\boc$ under which the left 
$\KK_G(X\T X)$-module structure on $\hat\KK_G(X\T X)$ corresponds to the left 
$\JJ_\boc$-module structure on $\hat\JJ_\boc$. If $\G\in L$, there is a unique
$\o'\in G\bsl X$ (the set of $G$-orbits in $X$) such that $\hat\c$ carries 
$R_{\o'}$ (see 4.3) onto a (simple) $\JJ_\boc$-submodule $M_\G$ of 
$\hat\JJ_\G$ whose isomorphism class is independent of $\G$; we say that this
is the special $\JJ_\boc$-module. The $\JJ_\boc$-module $M_\G$ admits a basis 
$\{e_{\G'};\G'\in L\}$ in which any $t_u$ (with $u\in\boc$) acts by a matrix 
with all entries in $\ZZ_{\ge0}$, namely the basis corresponding to the basis
$\{\bar V^{\o,\o'};\o\in G\bsl X\}$ of $R_{\o'}$.

\enddocument